**Note on the Merit Factors of Sequences**
**N.A. Carella, November, 2005.**

*Abstract*: A method for estimating the merit factors of sequences will be provided. The result is also effective in determining the nonexistence of certain infinite collections of cyclic difference sets and cyclic matrices and associated binary sequences.

## 1 Introduction

A quantitative measure of the goodness of a sequence of complex numbers, called the merit factor, involve the evaluation of norms and integrals. The current methods used to compute this measure are largely specific to the families of sequences under investigation. So far the merit factors of several infinite families of sequences have been determined. In some cases the evaluations are exact, and in other cases the evaluations are asymptotic, see [JN92], [BN02], [BN04], [JB05]. In this note a method for estimating the merit factors and determining the existence of families of sequences of bounded aperiodic autocorrelation will be provided. As an application, it will be shown that the number of Barker sequences is finite. The technique also has applications in the theory of cyclic difference sets and cyclic matrices.

## 2 Fundamental Materials

A selection of elementary facts applicable to the theory of sequences of complex numbers and related objects are recalled for the reader convenience. An extensive survey of the theory of sequences of complex numbers is given in [HT98]. Usually the emphasis is on sequences over finite sets such as $S = \{\ 0,\ 1\ \}$, $\{\ -1,\ 1\ \}$, $\{\ -1,\ 0,\ 1\ \}$, and the set roots of unity $\mu_m = \{\ 1,\ \xi,\ \ldots,\ \xi^{m-1}\ \}$.

*Autocorrelations and Merit Factor*

**Definition** 1.   The periodic autocorrelation and aperiodic autocorrelation of a sequence $s_0,\ s_1,\ \ldots,\ s_{n-1}$ of complex numbers are defined by

$$(1) \qquad \theta(t) = \sum_{k=0}^{n-1} s_k \bar{s}_{k+t}\ \text{ and }\ \rho(t) = \sum_{k=0}^{n-t-1} s_k \bar{s}_{k+t}$$

respectively. The autocorrelations satisfy many the relations, among these relations are
(i) $\theta(0) = \rho(0)$,        (ii) $\theta(t) = \rho(t) + \rho(n-t)$, for $t \not\equiv 0 \bmod n$,      (iii) $\rho(t) = \bar{\rho}(-t)$.

**Definition** 2.   The Fourier transform of a sequence $s_0,\ s_1,\ \ldots,\ s_{n-1}$ of complex numbers, and the power spectrum (its absolute value) are given by

$$(2) \qquad s(z) = \sum_{k=0}^{n-1} s_k z^k\ ,$$

and

$$(3) \qquad \left| s(z) \right|^2 = \sum_{k=0}^{n-1} s_k \bar{s}_k + \sum_{k=1}^{n-1} \rho(k)(z^k + \bar{z}^k)\ ,$$

respectively, where the pair $z = e^{i2\pi x}$, $\bar{z} \in T = \mathbb{R}/\mathbb{Z}$ are complex conjugates. The spectrum $\left|\ s(z)\ \right|^2$ is a continuous function on the unit sphere $T$ and its constant term $\rho(0) = \left| s(0) \right|^2 = s_0 \bar{s}_0 + \cdots + s_{n-1} \bar{s}_{n-1}$ coincides with the $l_2$ norm of the sequence.

**Definition** 3.   The merit factor of a sequence $s_0,\ s_1,\ \ldots,\ s_{n-1}$ of complex numbers is defined by

$$(4) \qquad F = \rho^2(0)\left( 2\sum_{t=1}^{n-1} \rho^2(t) \right)^{-1}.$$

This formula expresses the merit factor in terms of the values of the aperiodic autocorrelation of the sequence. The spectrum function of the sequence facilitates another formulation of the merit factor. This emerges from the equality

$$(5) \qquad 2\sum_{t=1}^{n-1} \rho^2(t) = \int_0^1 \left( \left| s(\alpha) \right|^2 - \left| s(0) \right|^2 \right)^2 d\alpha = \int_0^1 \left| s(\alpha) \right|^4 d\alpha - \left| s(0) \right|^4.$$

Combining this into the discrete formula (4) gives an analytical formula

$$(6) \qquad F = \left\| s(z) \right\|_2^4 \left( \left\| s(z) \right\|_4^4 - \left\| s(z) \right\|_2^4 \right)^{-1}.$$

This version gives an intuitive interpretation of the merit factor: It measures the deviation of the $L_4$ norm of $|s(z)|^2$ from a flat spectrum of magnitude $|s(0)|^4 = \left\| s(z) \right\|_2^4$, which is the mean value of $|s(z)|^2$ over the unit sphere $T = \mathbb{R}/\mathbb{Z}$. An ideal sequence has a large and increasing merit factor as a function of $n$.

A survey of the merit factors of many important binary sequences is given in [JB05]. Extensive numerical data on binary sequences appear in [KR05], and computational techniques are given in [JN92], [BN00], [BN04], data on complex sequences are given in [BN05], etc. And a generalization to multivariable sequences appear in [GR05], and [PR05].

A central problem in the theory of sequences is the determination of sequences in the binary sequences space $\mathcal{B}_n = \{ \, s = (s_0, s_1, \ldots, s_{n-1}) : s_i = \pm 1 \, \}$ of size $2^n$, ($q^n$ in the $q$-adic case) that maximize (4) and (6). Equivalently, it calls for the minimization of the sums of squares $\rho_1^2 + \cdots + \rho_{n-1}^2$ (sidelobe components) of sequences in $\mathcal{B}_n$.

## 3 Characterization By Autocorrelations

Let $S \subset \mathbb{C}$ be a subset of complex numbers, and let $S^n = \{ \, s = (s_0, s_1, \ldots, s_{n-1}) : s_i \in S \, \} \subset \mathbb{C}^n$ be the corresponding sequences space.

The process of matching sequence to aperiodic autocorrelation investigated here is depicted in a diagram:

$$(7) \qquad \begin{array}{ccccc} \mathbb{C}^n & \xrightarrow{\ \mathfrak{I}_n\ } & \mathbb{C}[z] & \xrightarrow{\ \ell_2\ } & \mathbb{C}[z^{-1}] \\[4pt] & \ddots & & & \downarrow \iota \\[4pt] & & \mathbb{C}^n & \xrightarrow{\ \mathfrak{I}_{\pm n}\ } & \mathbb{C}[z^{-1}] \end{array}$$

In the upper branch, the Fourier transform $\mathfrak{I}_n$ maps a sequence $s = (s_0, s_1, \ldots, s_{n-1}) \in \mathbb{C}$ to its image polynomial $f(z) \in \mathbb{C}[z]$, and the $\ell_2$ norm maps a polynomial $f(z) \in \mathbb{C}[z]$ to its absolute value $f(z)\overline{f(z)} \in \mathbb{C}[z^{-1}]$. In the lower branch, the two-sided Fourier transform $\mathfrak{I}_{\pm n}$ maps an aperiodic autocorrelation vector $\rho = (\rho_0, \rho_1, \ldots, \rho_{n-1})$ to its image polynomial $\mathfrak{I}_{\pm n}(\rho) = \displaystyle\sum_{k=-n+1}^{n-1} \rho_k z^k \in \mathbb{C}[z^{-1}]$, and $\iota$ is the identity map. Thus $\ell_2 \circ \mathfrak{I}_n(s) = $ $\mathfrak{I}_{\pm n}(\rho)$. The difficulty in finding sequences from a specific family of sequences arises in computing the inverse images $s = (\ell_2 \circ \mathfrak{I}_n)^{-1}(\rho)$ in $\mathbb{C}^n$. Specifically, given an aperiodic autocorrelation vector $\rho = (\rho_0, \rho_1, \ldots, \rho_{n-1}) \in \mathbb{C}^n$, determine the sequences $s = (s_0, s_1, \ldots, s_{n-1})$ that correspond to its inverse image.

The inverse image $(\ell_2 \circ \mathfrak{I}_n)^{-1}(R) \subset S^n$ of a subset of aperiodic autocorrelation vectors $R \subset \mathbb{C}^n$ is clearly dependent on the followings:

(i) The cardinality $\#S$ of the subset $S$.

(ii) The maximal aperiodic autocorrelation, that is, $\max \{ \, | \, \rho(t) \, | : t \not\equiv 0 \bmod n \, \}$.

For example, the inverse image of the aperiodic autocorrelation vector $\rho = (\rho_0, \rho_1, \ldots, \rho_{n-1}) = (n, |n-1|, |n-2|, \ldots, |n-t|, \ldots, 1)$ in the binary sequences space $\mathcal{B}_n$ contains the all ones sequence $s_0, s_1, \ldots, s_{n-1} = 1, 1, \ldots, 1$, the alternating sequence $t_0, t_1, \ldots, t_{n-1} = 1, -1, 1, -1, \ldots, 1$, and possibly other sequences.

The aperiodic autocorrelation vectors of the all ones sequence and the alternating sequence maximizes the sum of squares $\rho_1^2 + \cdots + \rho_{n-1}^2$, so minimal merit factor of binary sequences over $\{ -1, 1 \}$ is $F_n = 3n^2 / (2n^3 - 3n^2 + n)$, which has the $L_4$ norm $\| s(z) \|_4^4 = n(2n+1)/3$. In contrast, the largest merit factor possible is unknown. The estimated maximal merit factor of any family of sequences has been changed a few times after discovery of new families of sequences with higher merit factors than previous cases, see [BN04]. The determination of $\limsup F$ over the set of all binary sequences remains an open problem.

***Definition* 4.** A *p*-sequence is a sequence over a finite set $S$ with an aperiodic autocorrelation vector $\rho = (\rho_0, \rho_1, \ldots, \rho_{n-1})$ such that $\rho_0 \in \mathbb{N}$, and $| \, \rho_k \, | \leq (p-1)/2$ for $t \not\equiv 0 \bmod n$, where $p$ is a prime.

In particular, this gives a characterization of sequences in terms of the aperiodic autocorrelation sidelobe $(p-1)/2 = \max \{ \, | \, \rho_k \, | : t \not\equiv 0 \bmod n \, \}$. The sequence to aperiodic autocorrelation diagram of *p*-sequences has the following form:

$$(8) \qquad \begin{array}{ccccc} S^n & \xrightarrow{\;\mathfrak{I}_n\;} & \mathbf{F}_p[z] & \xrightarrow{\;\ell_2\;} & \mathbf{F}_p[z^{-1}] \\[2mm] & \diagdown & & & \downarrow \iota \\[2mm] & & \mathbb{R} \times \mathbf{F}_p^{n-1} & \xrightarrow{\;\mathfrak{I}_{\pm n}\;} & \mathbf{F}_p[z^{-1}] \end{array}$$

Extensive of works have been directed toward the investigation of *p*-sequences over $S = \{ -1, 1 \}$, with small $p = 2, 3$ and so on. A problem about *p*-sequences, for $p = 3$, with the aperiodic autocorrelation vectors $\rho = (\rho_0, \rho_1, \ldots, \rho_{n-1}) \in \mathbb{N} \times \mathbf{F}_3^{n-1}$ will be considered in later on.

***Golay Conjecture* 5.** The merit factors $F$ of binary sequences have an absolute maximal $F \leq F_{max}$ for all $n$. Equivalently, the norms $\left\| s(z) \right\|_4 \geq (1+c)\sqrt{n}$ for all $n$, and $c > 0$ constant.

## 4 Estimates of Exponential Sums

The results in this section are of independent interest and are stated in the original notations or standard notations.

Let $a_0, a_1, \ldots, a_{N-1}$ be a list of complex numbers. The interest here is in finding an effective estimate of the exponential sum

$$(9) \qquad\qquad S(x) = \sum_{n=0}^{N-1} a_n z^n \, ,$$

where $z = e^{i2\pi x}$, and $x \in \mathbb{R}$. The large sieve inequality gives an average estimate of the sum $S(x)$ at $R$ points.

***Proposition* 6.** Let $x_0, x_1, \ldots, x_{R-1}$ be distinct real numbers in the unit interval $[0, 1]$, and let $\delta = \min \| x_i - x_j \|$, $i \neq j$. Then

$$(10) \qquad\qquad \sum_{0 \leq r < R} \left| S(x_r) \right|^2 \leq (N + \delta^{-1}) \sum_{0 \leq n < N} \left| a_n \right|^2 \, .$$

The large sieve inequality, which is essentially a generalization of Parseval formula

$$(11) \qquad\qquad \sum_{0 \leq n < N} \left| \hat{f}(x_n) \right|^2 = N \sum_{0 \leq n < N} \left| f_n \right|^2 \, ,$$

links analysis in discrete space to analysis in continuous space. Other versions of the large sieve inequality are also known, see [MV73] and related literature. This quickly leads to the average order of the sum

$$(12) \qquad\qquad \sum_{0 \leq n < N} \left| \hat{f}(x_n) \right|^2 = N \sum_{0 \leq n < N} \left| f_n \right|^2 \, ,$$

where $\hat{f}(x_n) = S(x_n)$.

The estimate below can be viewed as an extension of Weyl criterion to the exponential sum $S(x)$. The well known Weyl criterion is the same as (9) but has unit coefficients $a_0 =$ $a_1 = \cdots = a_{N-1} = 1$, see [KS74], [KV92], etc. for technical details and related exponential sums.

**Proposition 7.** Let $a_0, a_1, \ldots, a_{N-1}$ be a sequence of complex numbers of $l_2$ norm $\mid a_0 \mid^2 + \cdots + \mid a_{N-1} \mid^2 = N^{\nu/(1+\varepsilon)}$, where $\nu \geq 1$, and $\varepsilon > 0$, and suppose that $\mid c_n \mid \leq O(N^{\nu-1})$, $n \neq 0$. Then

$$(13) \qquad \sum_{n=0}^{N-1} a_n z^n = o(N^{\nu/2}),$$

where $z = e^{i2\pi x}$, and $x \in \mathbb{R}$ is an irrational number.

Proof: As usual write $c_n = \sum_{0 \leq k < N-n-1} a_k \overline{a}_{k+n}$, and use the assumptions on its absolute values to obtain

$$
(14) \qquad
\begin{aligned}
\left| \sum_{n=0}^{N-1} a_n z^n \right|^2 &= \sum_{n=0}^{N-1} a_n \overline{a}_n + \sum_{n=1}^{N-1} c_n (z^n + \overline{z}^n) \\
&= \sum_{n=0}^{N-1} a_n \overline{a}_n + \sum_{n=1}^{N-1} O(N^{\nu-1})(z^n + \overline{z}^n) \\
&= o(N^\nu) + O(N^{\nu-1}) \sum_{n=1}^{N-1} (z^n + \overline{z}^n).
\end{aligned}
$$

Applying Weyl criterion to the finite sum and simplifying return

$$(15) \qquad \left| \sum_{n=0}^{N-1} a_n z^n \right|^2 = o(N^\nu). \qquad\qquad \blacksquare$$

This result can also be proved using the large sieve inequality, see [SZ94, p. 79] and the references therein. This estimate is very much in line with the average order of the exponential sum $S(x)$.

**Proposition 8.** Let $\mid a_0 \mid^2 + \cdots + \mid a_{N-1} \mid^2 = cN$, $c > 0$ constant. If $\psi > 0$ is an increasing function of $N$, then the set of numbers $x \in [0, 1)$ such that

$$(16) \qquad \left| \sum_{n=0}^{N-1} a_n z^n \right|^2 \geq \psi(N) \sum_{n=0}^{N-1} a_n \overline{a}_n,$$

for all sufficiently large $N$ has measure zero.

The proof for the special case $a_0 = a_1 = \cdots = a_{N-1} = 1$ appears in [LG95, p. 42].

*Low Estimate of Exponential Sums*

An important result employed in estimating the lower estimates of certain exponentials sums is included here. The expression $f(x) = \Omega(g(x))$ denotes the condition $f(x) \geq |g(x)|$ for all sufficiently large $x > 0$.

**Theorem 9.** ([PE87]) Let $f$ be a periodic function of period 1. Assume that

(i) The function $f$ is piecewise continuous: there exist real numbers $0 = x_0 < \cdots < x_d = 1$ such that $|f(x) - f(y)| \leq c|x - y|$ whenever $x, y \in [x_i, x_{i+1}]$, and $c > 0$ constant.

(ii) The integral $\int_0^1 f(t)dt = 0$.

Let $b_n$ denotes the $n$th Fourier coefficient of the function $f$, let $\alpha$ be an irrational real number, and let $|b_n| >> n^{-1}$ for $n \equiv 0 \bmod m$, some $m$. Then

$$(17) \qquad \sum_{0 \leq n < N} f(n\alpha) = \Omega(\log^{1/2}(N)).$$

**Corollary 10.** Let $|a_n| \geq 1$ be constants, and let $\alpha$ be an irrational real number, then

$$(18) \qquad \sum_{0 \leq n < N} a_n e^{i2\pi n\alpha} = \Omega(\log^{1/2}(N)).$$

Proof: By linearity it is sufficient to consider the function $f(t) = a_m e^{i2\pi \beta mt}$, where $a_m$ is a constant, and $t \in [-1/2\beta, 1/2\beta)$, with $0 \neq \beta \in \mathbb{R}$ constant. The corresponding Fourier series is

$$(19) \qquad f(t) = \sum_{n=-\infty}^{\infty} b_n e^{i2\pi \beta t} = \sum_{n=-\infty}^{\infty} a_m \frac{\sin \pi(m-n)\beta}{\pi(m-n)\beta} e^{i2\pi \beta nt}.$$

Clearly, $|f(x) - f(y)| \leq |a_m| |e^{i2\pi \beta mx} - e^{i2\pi \beta my}| \leq 2\pi m |a_m| |x - y| = c|x - y|$, and its Fourier coefficients $|b_n| = |a_n| >> n^{-1}$ satisfy the conditions of Theorem 9. The claim follows from these. ∎

This is a weak estimate but it is effective in the application considered here. In summary, all these results show that the sum $S(x)$ follows a square root cancellation rule, that it,

$$(20) \qquad \left| \sum_{n=0}^{N-1} a_n z^n \right| = O\left( \sum_{n=0}^{N-1} a_n \overline{a}_n \right)^{1/2 + \varepsilon}$$

for almost all $x \in [0, 1)$, and $\varepsilon > 0$.

## 5 Estimates of the Integral

The integral of interest here is precisely the $L_4$ norm

$$(21) \qquad \left\| f(z) \right\|_4^4 = \int_0^1 \left| f(t) \right|^4 dt$$

of Fourier transform $f(z)$ of the sequence $a_0, a_1, \ldots, a_{N-1}$ on the unit disk $D = \{ z \in \mathbb{C} : |z| = 1 \}$. The evaluation of (21) and related norms are challenging problems. Several cases have been solved already, see [JN92], [BS90], [BN00], [BN01], [BN03], [BN04], etc for the case $|a_i| = \pm 1$, and and [GV05] for certain cases of complex sequences $|a_i| = \pm 1$. Two slightly different estimates of the $L_4$ norm of complex sequences characterized by the aperiodic autocorrelation $|c_n| \le c$ will be considered in this section.

The analytical tools, such as the circle method, utilized to estimate related but much more difficult integrals such as the $k$th moment of the zeta function, are well established, see [VA97], [IV85, p.129] and related literature.

By the monotonicity of the $L_p$ norms, it follows that

$$(22) \qquad \left\| f(z) \right\|_1 \le \left| f(0) \right| = \left\| f(z) \right\|_2 \le \left\| f(z) \right\|_4 \le \cdots \le \left\| f(z) \right\|_\infty .$$

Accordingly $c_0 = |f(0)|^2$ is a lower estimate of $\left\| f(z) \right\|_4$. On the other hand, there is the trivial upper estimate $\left\| f(z) \right\|_4^4 \le c_0^2 + 2 \sum_{0 \le n < N} |c_n|^2$. A nontrivial estimate of (21) attached to a sequence $a_0, a_1, \ldots, a_{N-1}$ of complex numbers of bounded aperiodic autocorrelation $|c_i| \le c$ is of the form

$$(23) \qquad \left\| f(z) \right\|_4^4 = \left| f(z) \right|^2 + E(N) ,$$

where the error/oscillation term $E(N) = o(N^2)$ or better.

*Estimate of the Integral*

Recall that if a function $f(t)$ is defined everywhere on a measurable set $I$, and $Z \subset I$ is a subset of measure zero then the integrals of the function on the sets $I$ and $J = I - Z$ are the same. This observation will also be used in the next result.

**Proposition 11.** Let $|f(0)|^2 = |a_0|^2 + \cdots + |a_{N-1}|^2$, and let $c = \max\{ |c_1|, \ldots, |c_{N-1}| \}$, $c > 0$ constant. Then

$$(24) \qquad \int_0^1 \left| f(\alpha) \right|^4 d\alpha = [ |f(0)|^2 + o(N) ]^2 .$$

Proof: Use the expanded spectrum function of the sequence $a_0, a_1, \ldots, a_{N-1}$ as in (3) and apply Proposition 6, Corollary 9 to obtain

$$(25) \qquad \left| f(\alpha) \right|^2 = \left| f(0) \right|^2 + \sum_{n=1}^{N-1} c_n (e^{i2\pi n\alpha} + e^{-i2\pi n\alpha}) = \left| f(0) \right|^2 + o(N),$$

where $\alpha$ is an irrational real number. Now replace it in the integral:

$$(26) \qquad \int_I \left| f(\alpha) \right|^4 d\alpha = \int_J \left| f(\alpha) \right|^4 d\alpha = \int_J [\, | f(0) |^2 + o(N)]^2 \, d\alpha = [\, | f(0) |^2 + o(N)]^2,$$

where $I = (0, 1)$ and $J = I - Z$ is the subset of irrational numbers modulo 1. ∎

The evaluations of the integral (21) by other means and Weyl criterion give essentially the same estimates.

**Proposition 12.** Let $f(z) = a_{N-1}z^{N-1} + \cdots + a_1 z + a_0 \in \mathbb{C}[z], |z| = 1$. Then

$$(27) \qquad \int_{|z|=1} \left| f(t) \right|^4 dt = [\, | f(0) |^2 + 2 \sum_{1 \le n < N} c_n \cos(2\pi n\xi) \,]^2.$$

Proof: The function $F(x) = \int_0^x \left| f(t) \right|^4 dt$ is continuous and differentiable on the unit interval $(0, 1)$, in fact its derivative is $F'(x) = | f(x) |^4$. Thus by *Mean Value Theorem*, it follows that

$$(28) \qquad \int_I \left| f(t) \right|^4 dt = \int_J \left| f(t) \right|^4 dt = \int_J \left| f(t) \overline{f(t)} \right|^2 dt = \left| f(\xi) \overline{f(\xi)} \right|^2,$$

for some irrational number $\xi \in J = I - Z$. ∎

Note that in applications of the Mean Value Theorem, the differentiability condition is required since there are continuous functions that are not differentiable anywhere on an interval.

This result is applicable to any sequence of complex numbers. Moreover, if the aperiodic autocorrelation satisfies $| c_k | \le c$, then the term $o(N) = 2 \sum_{1 \le n < N} c_n \cos(2\pi n\xi) > 0$ for all $N$.

*Quasi Monte Carlo Estimate of the Integral*
There are other methods of approximating the integral (21) such as the $L_2$ triangle inequality. Another possible method is given in [VA97, p. 14]. A numerical estimate of (21) is considered here. The numerical analysis has threads in common with the previous one.

The quasi Monte Carlo estimate of the integral is by means of the Koksma-Hlawka inequality

$$(29) \qquad \left| \frac{1}{N} \sum_{n=0}^{N-1} f(x_n) - \int_0^1 f(t)\, dt \right| \leq V(f) D_N(\omega),$$

where $V(f)$ is the total variation of the function $f(t)$ over the unit interval $[0, 1)$ , and $D_N(\omega)$ is the discrepancy of the nodes $x_0, x_1, \ldots, x_{N-1}$ in $[0, 1)$. Other related forms of this inequality can also be used, see [KS74, p. 142].

Moreover, the fact that the integral of Fourier series of a function is easily converted into an effective numerical integration facilitates the computation of the estimate of (21).

The discrepancy of a sequence of real numbers $x_0, x_1, \ldots, x_{N-1}$ in $[0, 1)$ is a metric in the range $1/N$ to 1, that measures its deviation from an ideal uniform sequence; an exact formula for computing the discrepancy is given in [KS74, p. 91]. Also, there are methods for selecting suitable sequences. One of these methods is based on Diophantine approximation theory.

Algebraic irrational numbers $\omega \in \mathbb{R}$ have type $\eta \geq 1$, see the literature for precise definitions and other technical details.

***Proposition* 13.** ([KS74, p. 123]) Let $\omega$ be an irrational real number of type $\eta$. Then for every $\varepsilon > 0$, the discrepancy $D_N(\omega)$ of the sequence of fractional part $x_n = \{\omega n\}$, $n = 1, 2, \ldots, N$, satisfies $D_N(\omega) = O(N^{-1/\eta + \varepsilon})$.

***Proposition* 14.** ([KS74, p. 143]) Let $x_0, x_1, \ldots, x_{N-1}$ be a sequence of real numbers of discrepancy $D_N(\omega)$. Then

$$(30) \qquad \sum_{n=0}^{N-1} e^{i 2\pi x_n} \leq 2N D_N(\omega).$$

This is a fundamental result in the theory of uniform distribution of sequences related to Weyl criterion. This and the associated exponential sums are also discussed in [NR92], [KV92, p. 176], etc.

Without loss in generality, put $|f(0)|^2 = |a_0|^2 + \cdots + |a_{N-1}|^2 = N$.

**Proposition 15.** Let $| f(0) |^2 = N$, and let $c_1, \ldots, c_{N-1}$ be real numbers such that $| c_n | \leq c$ for $n \not\equiv 0 \bmod N$, where $c > 0$ is a constant. Then

$$(31) \qquad \int_0^1 \big| f(x) \big|^4 \, dx = N^2 + O(N^{1+\varepsilon}) \, ,$$

for every $\varepsilon > 0$.

Proof: Since $| f(x) |^2 \leq N + 4c(N-1) \leq 5cN$, (by hypothesis), the total variation of the function $g(x) = | f(x) |^4$ over $[0, 1)$ satisfies the inequality $V(g) \leq (5cN)^2$. Now choose the nodes in $[0, 1)$ of discrepancy $D_N(\omega) = O(N^{-1+\varepsilon})$, see Proposition 13. Then the error term in the approximation satisfies $R_N \leq V(f) D_N(\omega) \leq O(N^{1+\varepsilon})$.

Now rewrite (21) as a quadrature formula:

$$(32) \qquad \int_0^1 \big| f(x) \big|^4 \, dx = \frac{1}{N} \sum_{n=0}^{N-1} \big| f(x_n) \big|^4 + R_N \, .$$

Replacing the expanded spectrum (3) and using the fact that the $c_n$ are real numbers return

$$
\begin{aligned}
(33) \qquad \frac{1}{N} \sum_{m=0}^{N-1} \big| f(x_m) \big|^4 &= \frac{1}{N} \sum_{m=0}^{N-1} \left| N + \sum_{n=1}^{N-1} c_n (z^n + \bar{z}^n) \right|^2 \\
&= \frac{1}{N} \sum_{m=0}^{N-1} \left( N^2 + 2N \sum_{n=1}^{N-1} c_n (z^n + \bar{z}^n) + \left( \sum_{n=1}^{N-1} c_n (z^n + \bar{z}^n) \right)^2 \right) \\
&= N^2 + 2 \sum_{m=0}^{N-1} \sum_{k=1}^{N-1} c_n (z^n + \bar{z}^n) + \frac{1}{N} \sum_{m=0}^{N-1} \left( \sum_{n=1}^{N-1} c_n (z^n + \bar{z}^n) \right)^2 .
\end{aligned}
$$

The first multiple finite sum reduces to

$$(34) \qquad 2 \sum_{m=0}^{N-1} \sum_{n=1}^{N-1} c_n (z^n + \bar{z}^n) = 2 \sum_{n=1}^{N-1} c_n \sum_{m=0}^{N-1} (e^{i2\pi n x_m} + e^{-i2\pi n x_m}) = O(N^{1+\varepsilon}) \, .$$

Here the inner sum has a magnitude of $O(N^\varepsilon)$, see Proposition 13, and outer sum has a magnitude of $O(N)$. Similarly, the second multiple sum reduces to

(35)
$$\frac{1}{N}\sum_{m=0}^{N-1}\left(\sum_{n=1}^{N-1}c_n(z^n+\bar{z}^n)\right)^2=\frac{1}{N}\sum_{n=1}^{N-1}c_n(e^{i2\pi nx_m}+e^{-i2\pi nx_m})\sum_{n=1}^{N-1}c_n\sum_{m=0}^{N-1}(e^{i2\pi nx_m}+e^{-i2\pi nx_m})$$
$$=O(1)\sum_{n=1}^{N-1}c_n\sum_{m=0}^{N-1}(e^{i2\pi nx_m}+e^{-i2\pi nx_m})=O(N^{1+\varepsilon}).$$

For sufficiently large $N$ (or asymptotically), the error term in (32) is absorbed in the term $O(N^{1+\varepsilon})$ or conversely. Therefore the estimate of the integral reduces to

(36)
$$\int_0^1\left|f(x)\right|^4dx=\frac{1}{N}\sum_{m=0}^{N-1}\left|f(x_m)\right|^4+R_N=N^2+O(N^{1+\varepsilon}).\qquad\blacksquare$$

It should be noted that all these estimates (24), (31), etc of (21) are almost identical, (have the same order of magnitude) and any one of these can be used in the applications given here.

## 6 Applications to Barker Sequences

This section considers the first application of the new estimates of the integral (21).

**Definition 16.** A binary sequence $s_0, s_1, \ldots, s_{n-1}$ over the set $\{-1, 1\}$ is called a Barker sequence if the aperiodic autocorrelation satisfies $\rho(0) = n$, and $|\rho(t)| \leq 1$ for $0 < t < n$.

Barker sequences are characterized by the aperiodic autocorrelation $|c_k| \leq 1$. Similar collections of sequences are characterized by the aperiodic autocorrelations $|c_k| \leq c$, $c > 0$ constants. The entire collection of known Barker sequences are listed here, see [TN65, p. 394].

The initial tail of the collection $|\rho(t)| \leq 1$.

| | |
|---|---|
| $n = 2$ | $s_k = 1, 1$ |
| $n = 3$ | $s_k = 1, 1, -1$ |
| $n = 4$ | $s_k = 1, 1, 1, -1$ and $1, 1, -1, 1$ |
| $n = 5$ | $s_k = 1, 1, 1, -1, 1$ |
| $n = 7$ | $s_k = 1, 1, 1, -1, -1, 1, 1-1$ |
| $n = 11$ | $s_k = 1, 1, 1, -1, -1, -1, 1-1, -1, 1, -1$ |
| $n = 13$ | $s_k = 1, 1, 1, 1, 1, -1, -1, 1, 1, -1, 1-1, 1$ |

Three other equivalent sequences are given by the transformations

i. $t_k = -s_k$,        ii. $t_k = (-1)^k s_k$,        iii. $t_k = (-1)^{k+1}s_k$.

**Proposition 17.** The aperiodic autocorrelation $c_k$, $0 < k < n$, of a Barker sequence of length $n > 0$ satisfies the relation

(37)
$$c_k = \begin{cases} 0 & \text{if } n-k \text{ is even,} \\ \pm 1 & \text{if } n-k \text{ is odd.} \end{cases}$$

This fact stems from the relations $c_k + c_{n-k} \equiv n \bmod 4$, $c_k \equiv n-k \bmod 4$, and $|c_k| \le 1$. For more details, see [TN61], [BN02, p.101].

Quite often a Barker sequence is given in polynomial form (its Fourier series), for example, the first few are expressed as $f_2(z) = z + 1$, $f_3(z) = z^2 + z - 1$, ,$f_4(z) = z^3 + z^2 + z - 1$, ... . The sequence of Barker polynomials minimize the $L_4$ norm on the unit disk $D = \{ z \in \mathbb{C} : |z| = 1 \}$.

**Proposition** 18. The $L_4$ norm of a Barker sequence of length $n > 0$ satisfies the relation

(38)
$$\left\| f(z) \right\|_4^4 = \begin{cases} n^2 + n & \text{if } n \text{ is even,} \\ n^2 + n - 1 & \text{if } n \text{ is odd.} \end{cases}$$

Proof: Use $\left\| f(z) \right\|_4^4 = c_0^2 + 2(c_1^2 + \cdots + c_{n-1}^2)$, $c_0 = n$, and Proposition 17. ∎

The merit factors of these sequences range from 1 to $169/12 = 14.083$. The precise equation of the merit factor of Barker sequences is stated here.

**Proposition** 19. The merit factor of a Barker sequence of length $n > 0$ is given by

(39)
$$F_{\mathrm{B}} = \begin{cases} n & \text{if } n \text{ is even,} \\ n^2/(n-1) & \text{if } n \text{ is odd.} \end{cases}$$

Proof: It follows from the discrete formula (4), and Proposition 18. ∎

**Conjecture** 20. There is no Barker sequence of length $n > 13$.

The current best results state that there are neither Barker sequences of odd length $n > 13$, see [TN61], nor even length $n < 4 \cdot 10^{12}$, see [ST02, p. 73]. The method utilized here gives a unified proof for all large $n$ independent of the parity of $n$ and other arithmetic properties.

**Theorem** 21. The number of Barker sequence of length $n \ge 1$ is finite.

Proof: In light of Proportion 11 or 15, and the sequence parameters $|f(0)|^2 = c_0 = n$, and $|c_k| \le 1$ for $0 < k < n$, it follows that the integral

(40)
$$\left\| f(z) \right\|_4^4 = \int_I \left| f(t) \right|^4 dt = \left[ n + o(n) \right]^2 .$$

Thus the merit factor of any such sequence satisfies the relation

(41)
$$F_B = \frac{n^2}{\left\| f(z) \right\|_4^4 - n^2} = \frac{n^2}{\left[ n + o(n) \right]^2 - n^2} = o(n)$$

But by Proposition 19, the merit factor also satisfies

(42)
$$\frac{n^2}{\left[ n + o(n) \right]^2 - n^2} = \begin{cases} n & \text{if } n \text{ is even,} \\ n^2/(n-1) & \text{if } n \text{ is odd.} \end{cases}$$

Here the error term $E(n) = \left[ n + o(n) \right]^2 - n^2$ grows at least as fast as $\left[ n + \Omega(\log(n)^{1/2}) \right]^2 - n^2$, see Corollary 10, but it is of order $o(n^2)$ as $n$ increases. Therefore in order to satisfy both constraints the order $n$ of the sequence must satisfy $n \leq n_0$, where $n_0 > 0$ is a constant. ∎

It is quite possible that by computing a more accurate estimate of the integral, a concrete upper limit on the order of the form $n < 10^K = n_0$ could be obtained. Then the Schmidt test, see [ST02, p. 71], could be used to settle the remaining finite number of cases.

### 7 Menon Cyclic Difference Sets and Cyclic Hadamard Matrices

For odd $v > 6$, the cyclic difference sets corresponding to Barker sequences are

$$(v, k, m) = (7, 4, 2), (11, 5, 2) \text{ and } (13, 9, 6).$$

For even $v > 13$, a Barker sequence exists if and only if there is a cyclic Menon difference set $(v, k, m) = (4u^2, 2u^2 \pm u^2, u^2 \pm u^2, u^2)$, $u \geq 2$, see [BT71, p. 97].

***Theorem 22.*** The number of cyclic Menon difference sets $(4u^2, 2u^2 \pm u^2, u^2 \pm u^2, u^2)$, $u > 0$ is finite.

A cyclic Hadamard matrix $H = \langle a_0, a_1, \ldots, a_{n-1} \rangle$ with $a_i \in \{ -1, 1 \}$, which is characterized by $HH^T = I_n$, is generated by any one of its rows. Only one such matrix is currently known, namely, $H = \langle 1, 1, 1, -1 \rangle$.

***Conjecture 23.*** There is no cyclic Hadamard matrix of order $n = 4u^2 > 4$.

The current best result states that there is no cyclic Hadamard matrix of order $n < 4 \cdot 10^{11}$, with possible exceptions of $u = 165, 11715,$ and $82005$, see [ST02, p. 73].

A cyclic Hadamard matrix of order $n = 4u^2$ exists if and only if there is a cyclic Menon difference set $(v, k, m) = (4u^2, 2u^2 \pm u^2, u^2 \pm u^2, u^2)$, $u \geq 2$, see [BT71, p. 97]. On the other hand, it is stated in [JL99, p.423] that there is no proof of this claim, see also [PT95, p. 152] for a similar discussion. A well known result of relevant to this discussion is restated here:

**Proposition 24.** ([JL99]) A periodic binary sequence with period $v$, $k$ entries +1 per period and two-level autocorrelation function (with all nontrivial autocorrelation coefficients equal to γ) is equivalent to a cyclic $(v, k, m)$-difference set, where $\gamma = v - 4(k - m)$.

For example, a binary sequence $a_0, a_1, \ldots, a_{n-1}$ generated by cyclic Menon difference set has a two levels periodic autocorrelation $\theta(0) = v = 4u^2$, and $\theta(t) = v - 4(k - m) = 0$ for $t \not\equiv 0 \bmod v$.

It appears that the argument in [JL99] is asking for a proof of the following: Does the two levels periodic autocorrelation

$$(43) \qquad \theta(t) = \sum_{i=0}^{n-1} a_i a_{i+t} = \begin{cases} v = 4u^2 & \text{if } t = 0, \\ 0 & \text{if } t \neq 0, \end{cases}$$

imply the four-level aperiodic autocorrelation function

$$(44) \qquad \rho(t) = \sum_{i=0}^{n-1-t} a_i a_{i+t} = \begin{cases} v = 4u^2 & \text{if } t = 0, \\ 0, \pm 1 & \text{if } t \neq 0. \end{cases}$$

Observe that for $v = \text{prime} \neq 4u^2$ this condition does not hold, for instance, the Legendre sequence gives a counterexample.

**Theorem 25.** Assuming the above, the number of cyclic Hadamard matrix is finite.

A sequence is called *perfect* if $\theta(0) \neq 0$, and $\theta(t) = 0$ for $t \not\equiv 0 \bmod n$.

**Corollary 26.** The number of perfect sequences over $\{ -1, 1 \}$ is finite.
In contrast to the binary case, perfect sequences over the complex numbers are abundant.

**Proposition 27.** (Turyn 1968) For any odd $n > 1$ there exits a perfect sequence over the set root of unity $\mu_n = \{ 1, \xi, \ldots, \xi^{n-1} \}$.

**8 Sequences with Small Aperiodic Autocorrelations**

The characterization of the merit factors of a family of sequences by the aperiodic autocorrelation is natural. Also note that the merit factors of a family of sequences can not be characterized by the periodic autocorrelation.

The analysis employed in the special case of binary sequences of aperiodic autocorrelation $| \rho(t) | \leq 1$ covered so far shed light on the problem of finding sequences of length $n$ over the set $\{ -1, 1 \}$, which have small aperiodic autocorrelation $| \rho(t) | \leq c$, for $t \not\equiv 0 \bmod n$, where $c > 0$ is a constant.

The next best family of binary sequences after the collection of Barker sequences characterized by $| \rho(t) | \leq 1$ for $t \not\equiv 0 \bmod n$ is the collection of binary sequences characterized by $| \rho(t) | \leq 2$ for $t \not\equiv 0 \bmod n$. This too is a finite collection, but a larger collection that includes the Barker sequences as a subset. A few of these are listed here.

The initial tail of the collection $| \rho(t) | \leq 2$.

$n = 2$            $s_k = 1, 1$

$n = 3$            $s_k = 1, 1, 1$; and every binary sequence of length 3.

$n = 4$            $s_k = 1, 1, 1, 1$; and every binary sequence of length 4, but not $t_k = 1, 1, 1, 1$; $t_k = -1, -1, -1, -1$.

$n = 5$            $s_k = 1, 1, 1, 1, -1$; and every binary sequence of length 5 but not

$t_k = 1, 1, -1, -1, 1$;        $t_k = 1, -1, 1, -1, 1$;        $t_k = 1, -1, -1, 1, 1$;

$t_k = -1, 1, 1, -1, -1$;       $t_k = -1, -1, 1, 1, -1$;       $t_k = -1, -1, -1, -1, -1$.

***Proposition* 28.**     Suppose that a sequence $s_0, s_1, \ldots, s_{n-1}$ of complex numbers has bounded aperiodic autocorrelation $| \rho(t) | \leq c$, for $t \not\equiv 0 \bmod n$, where $c > 0$ is a small constant. Then its merit factor satisfies the inequality $F \geq n^2/[2(n-1)c] \geq n/(2c^2)$.

***Theorem* 29.**     The set of sequences $s_0, s_1, \ldots, s_{n-1}$ over a finite set $S \subset \mathbb{C}$ of complex numbers, and bounded aperiodic autocorrelation $| \rho(t) | \leq c$, for $0 < t < n$, where $c > 0$ is a constant, is finite.

This generalizes Theorem 21. The existence of an infinite collection sequences of bounded aperiodic autocorrelations $| \rho(t) | \leq c$, for $t \not\equiv 0 \bmod n$, implies that the merit factors of orders $o(n)$, using Proposition 11. But this contradicts Proposition 28. Thus, for each fixed $c > 0$, the collection is finite.

Not surprisingly, all the families of sequences over finite sets of complex numbers of small autocorrelation and large merit factors reported in the literature are finite, see [TN61], [TN74], [BN05], etc.

**9 Sequences of Large Aperiodic Autocorrelations**

Sequences of small periodic autocorrelations $\theta(t)$ can have large aperiodic autocorrelations $\rho(t)$, but not conversely. This follows from $\theta(t) = \rho(t) + \rho(n - t)$, for $t \not\equiv 0 \bmod n$, For example, the sequences generated by cyclic difference sets of small periodic autocorrelations can have large aperiodic autocorrelations.

The important infinite collection of character sequences $\{ \chi(t) : 0 \leq t < n \}$, where $\chi$ is a multiplicative character modulo $n$, do not have bounded aperiodic autocorrelations. In fact, the estimate

$$(45) \qquad \rho(t) = \sum_{k=1}^{n-t-1} \chi(k)\chi(k+t) = O(n^{1/2}\log(n)),$$

is close to the best possible, see [MV77]. As a result, the sums of squares are fairly large and of the form

$$(46) \qquad 2\sum_{t=1}^{n-1} \rho^2(t) = an^2 + o(n^2),$$

where $a > 0$ is a small parameter depending on the sequence, see the literature.

Assuming sequences over a finite set $S \subset \mathbb{C}$ of complex numbers, unbounded aperiodic autocorrelations $|\rho(t)|$, $0 < t < n$, account for the small merit factors of the sequences. This is clear from (4), which have large sums of squares. Conversely, bounded aperiodic autocorrelations $|\rho(t)|$, $0 < t < n$, account for the large merit factors of the sequences.

In contrast, sequences over an unbounded set $S \subset \mathbb{C}$ of complex numbers, can have unbounded merit factors. For example, the well known sequence, expressed as a polynomial,

$$(47) \qquad f(z) = \sum_{k=0}^{n-1} e^{ik(k+1)\pi/n} z^k,$$

has unbounded merit factor $F = o(n)$. This follows from the norm

$$(48) \qquad \left\| f(z) \right\|_4^4 = 2\sum_{t=0}^{n-1} \rho^2(t) = n^2 + o(n^2),$$

and (6) see [LD66].

**Distribution Of Merit Factors**

It is not difficult to show that the merit factors of all but balanced or nearly balanced binary sequences are bounded. Further the Central limit theorem/Law of large number seems to imply the existence of binary sequences of unbounded merit factors but extremely scarce. The distribution of the merit factor appears to be skewed symmetric of mean $\mu = 1$, such that the unbalanced binary sequences of bounded merit factors occur in the short left tail, and the main portion of the graph. The remaining sequences of unbounded merit factors occur on the long right tail.

***Open Problem* 30.** Is the collection of sequences $s_0$, $s_1$, …, $s_{n-1}$ over { −1, 1 } of bounded aperiodic autocorrelations $| \rho(t) | \leq c\log(n)$ for $t \not\equiv 0 \bmod n$, where $c > 0$ is a constant, infinite?

***Open Problem* 31.** Show that it is sufficient to prove Golay merit factor conjecture for balanced or nearly balanced binary sequences.